\pdfoutput=1

\input mypapers.mac
\input auxi.mac

\def\content{{\rm Spec}}

\def\content{{\rm Spec}}

\documentclass{hebtex}

\begin{document}

\norevhebuse 
\nolabelnote

\twelverom   
\automathbaselineskip 

\forcecolortrue 


\begintitle
Infinite sets of Mutually Coprime Locally Square-free 
Elements, Inside the Range of Polynomials with Values in Principal 
Ideal Domains 
\endtitle

\author{Michael Bensimhoun}{}{}{\date, Jerusalem}

\beginabstract

Let §\calR§ be an infinite principal ideal domain and §\calK = {\rm quot}(\calR)§. 
Assume that §P_1,\ldots P_n\in \calK[X]§ are polynomials which take §\calR§ 
to §\calR§, and §P§ is their product. If the §P_i§ satisfy
necessary conditions, there exists an infinite set §S§ 
such that the elements §P(m)§ are mutually coprime as §m§ varies in §S§,
and §P_i(m)§ is coprime to §P_j(m)§ for every §i\ne j§.
We prove that, 
in addition to the above property, if §P_i§ has no multiple roots
whenever §i§ belongs to some subset of §\set{1,\ldots n}§, 
§S§ can be constructed in such a way
that §P_i(m)§ is either a unit, or is divisible by some prime §p§, but
not by §p^2§. 
This result is the basis of a conjecture formulated at the end of this 
article, according to which one can extract infinitely many square-free 
elements from the value set of §P§, provided §P§ has no multiple~root. 
In the course of this article, the notion of \emph{totally primitive} polynomial
is introduced and elaborated in order to provide a convenient framework 
for those questions, or other questions concerning the value set of polynomials.  
This framework makes
possible more enlightening, and seemingly more 
general, statements of the famous Bunyakovsky and Schinzel conjectures. 
\endabstract
\medbreak

\keywords Polynomial, coprime, integer-valued, Bunyakovsky conjecture, Schinzel hypothesis


\begincontent
\item{\sref{intro}.} Introduction 
\item{\sref{def}.} Definitions and setup 
\item{\sref{thm}.} Main theorem and corollaries 
\item{\sref{proof}.} Proof of theorem \proclref{main} 
\item{\sref{necessary_cond}.} Note on prime moduli annihilating a polynomial 
\item{\sref{conjectures}.} Bunyakovsky and Schinzel conjectures restated and extended 
\item{\sref{PrincConj}.} Related conjectures
\endcontent



\section intro introduction

The aim of this article is two-fold:
on the one hand, it provides theorems on the coprimality of subsets 
of the value set of polynomials, and perhaps more importantly, on 
their ``locally square-free'' behaviors. 
These results, by themselves, could be useful in several circumstances.  
It should be observed immediately that nothing here is essentially new 
(and to be honest, not elementary): in its most restricted form, the main theorem 
in this article was known to Schinzel (and probably long before him), 
while Sury, in the integer case, 
already provided most of the ingredients needed for the proof (\cite{Su}). 
In some sense, this was rather a matter of putting everything
together. 

On the other hand, this article aims at providing a coherent framework 
for questions related to the value set of polynomials. 
Building upon the notion of \emph{total primitivity},  
few theorems allow manipulating easily the basic concepts.  
This framework leads eventually to more enlightening and extended exposition of the 
Bunyakovski and Schinzel conjectures. 

There was no reason to limit the extent of these results to the integer case.
In fact, one is naturally led to consider polynomials over the quotient ring of 
an integral domain §\calR§, which take §\calR§ to itself, 
by the theory of integer-valued polynomials.

Those are polynomials §P§ of §\bbQ[X]§ such that §P(\bbZ)\incl \bbZ§.
Every polynomial with integral coefficients is of course integer-valued, 
but integer-valued polynomials may belong to §\bbQ[X]\diff \bbZ[X]§,
as, for example, the so-called ``binomial polynomials'' 
$$
{X \binomial k} = {X(X-1)\ldots (X-k+1)\over k!}.
$$

Unfortunately, in order to generalize to §\calR§ the results in the integer case, 
most of the arithmetic properties of §\bbZ§ are needed,  
especially the prime factorization and the Chinese remainder theorems. 
This limits seriously the choice of §\calR§. Overall, principal ideal
domains seem to be the best candidate for an extension of the integer case theory. 

At the end of this article, we propose a (possibly new) conjecture, in 
relation with the Bunyakovski and Schinzel conjectures, and with 
the results exposed in this article. 
We think it is much easiest than the former, and though, could 
provide an impetus for the development of far reaching methods.   


\section def Definitions and setup

In this article, we assume that §\calR§ is an \emph{infinite}\footn{In other words, 
it is not a finite field.} principal ideal domain, whose quotient field is §\calK§,
and denote by §\calU§ the group of units of §\calR§. 
The ``modulo'' notation will be used in two different, yet not confusing, ways:
if §a,b\in \calR§ and §b = au§ with §u\in \calU§, we note §b = a\mod \calU§;
on the other hand, if §c\in \calR§ and §b-a = kc§ with §k\in \calR§, we note, as usual, §b = a \mod c§.

The set of prime elements, modulo §\calU§, dividing at least one element of a given set
§\set{k_1,k_2\ldots}\incl \calR§
will be noted §{\bf Spec}\set{k_1, k_2,\ldots}§.
So, §\content\, \calR§ is the set of prime elements of §\calR§ modulo §\calU§.

A polynomial §P\in \calK[X]§ is {\bf integrally valued}, or {\bf §\calR§-valued}, if §P(\calR)\incl \calR§.
We denote by §\bbI_\calR[X]§ the set of integrally valued polynomials with respect to §\calR§.
Thus, if §R = \bbZ§, ``integrally valued polynomial'' means ``integer valued polynomial''.

If §P\in \bbI_\calR[X]§,
§\calP(P)§ will denote the (finite) set of primes modulo §\calU§ that divide 
§P(m)§, for every §m\in \calR§.
Notice that §\calP(P) = \emptyset§ is equivalent to §\gcd P(\calR) = 1 \mod \calU§.  

A polynomial §P\in \bbI_\calR[X]§ is {\bf totally primitive} if §\calP(P) = \emptyset§.
Thus, a totally primitive \emph{constant} polynomial \emph{must} be a unit of §\calR§.

\remark
One may wonder how easy could it be to produce
non trivial integrally valued polynomials, 
in principal ideal domains other than 
§\bbZ§. A simple example is the following.

Let §p§ be a prime number, §\calR§ the localization of §\bbZ§ at §p§, and §\calK§ the quotient field of $\calR$. 
Of course, $\calR$ is a principal ideal domain, whose unique prime, up to a 
unit of $\calR$, is $p$.

Let §P = {1\over p} (X^p - X)§. Thus, §P§ belongs to §\calK[X]§ but not to §\calR[X]§.
Observe that the ``modulo p'' homomorphism $\bbZ \to \bbZ_{/p}$ extends canonically to a homomorphism 
$\calR\to \bbZ_{/p}$, hence §m^p - m\congr 0 \mod p§ for every $m\in \calR$.
This shows that  $P$ fulfills $$P(\calR)\incl \calR.$$ 

Another example: take $R=\bbF_q[T]$ and $\pi\in R$ an irreducible polynomial of degree $d$. 
Then $R/(\pi)$ is a finite field with $q^d$ elements. Hence $P={1\over \pi}(X^{q^d}-X)$ is integrally
valued. 

\proposition repr
Every polynomial §Q\in \calK[X]§ can be written in the form 
$$
Q(X) = {A\over B} P(X),
$$
where §A, B\in \calR§ are coprime, §P§ is totally primitive, 
and §A, B§ and §P§ are uniquely defined modulo §\calU§.

If §Q\in \bbI_\calR[X]§, then §B = 1\mod \calU§.
\endproclaim 

\proof 
let us write 
$$
Q = \alpha R(X)
$$
where §\alpha \in \calK[X]§ and §R\in \calR[X]§.
Denote §\gamma = \gcd R(\calR)§ (mod §\calU§), and 
$$
P(X) = {1\over \gamma}R.
$$ 
Then §P§ is integrally valued, and we set
$$
{A\over B} = \alpha\gamma,\quad \mtext{§A§ coprime to §B§.}
$$ 
This shows that §Q§ can be written in the contended form.
On the other hand, if 
$$
Q = q_1 P_1(X) = q_2 P_2(X),
$$
with §q_1,q_2\in \calK§, and §P_1, P_2§ totally primitive, 
then denoting §A/B = q_1/q_2§, with §A,B\in \calR§ coprime, there holds
$$
{A\over B}P_1(X) = P_2(X) \mod \calU.
$$ 
But §P_2(m)\in \calR§ for every §m\in \calR§, and we can find §m\in \calR§ such 
that §P_1(m)§ is coprime to any prime factor of §B§, as stated in 
\cor \proclref{almostNote} below.
This implies §B \in \calU§, and hence §A\in \calU§ since §P_2§ is totally primitive. 
\endproof

It is clear that any power of a totally primitive polynomial is totally primitive.
Moreover, if a product of integer valued polynomials §P_i§ is totally primitive,
then so are the §P_i§. 
 
Nevertheless, the product of totally primitive polynomials needs not be primitive; 
for example the product of §x§ and §x - 1§ is always divisible by §2§ for 
every §x\in \bbZ§, and the product of §x§, §x-1§ and §x-2§ is divisible by §3!§.
So, a totally primitive polynomial is not always the product of totally 
primitive irreducible factors.
Nevertheless, the following proposition follows easily from the remarks above.

\proposition unique_primitive_fact
A totally
primitive polynomial §P§ can be written uniquely modulo §\calU§, up to the order of the factors, 
in the form
$$
P(X) = {1\over \Gamma} P_1(X)^{\alpha_1} P_2(X)^{\alpha_2}\ldots P_n(X)^ {\alpha_n},
\quad \alpha_i \in \bbN^*,
$$
where §\Gamma \in \calR^\times§ and §P_1,P_2,\ldots § are totally primitive polynomials, 
irreducible in §\calK[X]§. Actually, denoting by §\Ptil§ the product of the §P_i^{\alpha_i}§,   
$$
\Gamma = \gcd \Ptil(\calR).
$$
\endproclaim

Joining the two propositions above, we conclude by the following result. 

\proposition primitive_facto 
A polynomial §Q\in \calK[X]§ can be written uniquely modulo §\calU§, up to the order of 
the factors, in the form
$$
Q(X) = {A\over B\Gamma} P_1(X)^{\alpha_1} P_2(X)^{\alpha_2}\ldots P_n(X)^ {\alpha_n} ,
\quad \alpha_i \in \bbN^*,
$$
with §A, B, \Gamma\in \calR§, §A§ coprime to §B§, §P_i§ totally primitive and irreducible, 
and 
$$\Gamma = \gcd \Qtil(\calR),$$
§\Qtil§ denoting the product of the §P_i^{\alpha_i}§.
\endproclaim

The above proposition shows that it makes sense to introduce de notions of 
{\bf totally primitive factorization} of 
§Q§, and of {\bf totally primitive (prime) factor} of~§Q§.  

It is now visible that any question concerning the value set of 
polynomials in §\calK[X]§ reduces to a question
about totally primitive polynomials.
They constitute the right framework, as far as these questions are concerned, 
and there is actually no point to proceed with general 
polynomials in §\calK[X]§.  
Henceforth, we shall consider only totally primitive polynomials in the sequel.


\section thm Main theorem and corollaries

The main theorem of this article, which will be proved in the next section, 
is an attempt to approach as much as possible the Bunyakovsky and Schinzel conjectures
in principal ideal domains.
Despite it is elementary, it may prove sufficiently powerful
in certain circumstances or practical cases, where these conjectures would be required.  

\theorem main
\listskipamount = 0pt
Let §P_1, P_2, \ldots P_n\in \calK[X]§ be non constant, totally primitive polynomials,
with respect to §\calR§.

Denote by §\calE§ a subset of 
§\set{1,2, \ldots n}§ (possibly a singleton) such that the polynomials §P_i(X)§ are pairwise 
coprime in §\calK[X]§, as §i§ belongs to §\calE§.
Further, denote by  §\calI§ the largest subset of  
§\calE§ such that §P_i§ 
has no multiple roots in an extension of §\calK§, 
as §i§ belongs to §\calI§. 
Finally, set 
$$
P = P_1P_2\cdots P_n.
$$
If §P§ is totally primitive, 
there exists an infinite set §S\subset \calR§ such that,
for every §m, m_1, m_2\in S§, with §m_1\ne m_2§,

\beginlist
\item{(i)}  
§P(m_1)§ is distinct from, and coprime to, §P(m_2)§;

\item{(ii)} §P_i(m_1)§ is distinct from (and coprime to) §P_i(m_2)§; 

\item{(iii)}
for every §i,j\in \calE§\/, with §i\ne j§, 
§P_i(m)§ is distinct from, and coprime to, §P_j(m)§. 

\item{(iv)}
for every §i \in \calI§,
§P_i(m)§ is either a unit, or 
is divisible by a certain prime §p = p(m,i)§, but not by~§p^2§.
Moreover, if §\calU§ is finite, and in particular if §\calR = \bbZ§, it can be 
supposed that §P_i(m)\nin \calU§; 

\item{(v)}
if  §\,\content\; \calR§ is finite, §P_i(m)§ can be furthermore assumed to belong to §\calU§, 
for every §i\in \calI§;

\item{(vi)}
if §\calR = \bbZ§,
§S§ can be chosen to be a subset of §\bbZ_+§, or of §\bbZ_-§.  
\endlist
\endproclaim

\remark
Assume that §P\in \bbZ[X]§ 
is a non constant polynomial
that has no roots modulo a prime §p\in \bbN§, like §X^2+1§ modulo §3§.
Let §\calR§ be the localization of §\bbZ§ at §p§.
Since the ``modulo p'' homomorphism extends canonically to §\calR§,
§P(m)§ cannot be a multiple of §p§ for every §m\in \calR§. 
So §P(m)§ is a unit for every §m\in \calR§.
This shows that in principal ideal domains, there may 
exist non constant integrally valued polynomials $P$ that take the whole domain to the units group,
while this is impossible in §\bbZ§. So, the two cases in (iv) are indeed possible.  
\medbreak

From (v), the following result, which can also be proved directly, is immediate.

\corollary infUnit
An infinite principal ideal domain which have only finitely many primes contains infinitely many units.
\endproclaim

The next corollary follows immediately from (ii) and the fact that §S§ is infinite.

\corollary almostNote
Assume that §E\incl \calR^*§ is finite.
For every §i§, the elements of §P_i(S)§ can be supposed coprime 
to the elements of §E§ in \thm \proclref{main}. 

Thus, a totally primitive and non constant polynomial has infinitely
many values coprime to the elements of any given finite set of numbers. 
\endproclaim

The next corollary allows manipulating the totally primitive factors of a polynomial.

\corollary FactManip 
Let §P\in \calK[X]§ be a non constant totally primitive polynomial, and §P_0§ be a totally primitive
prime factor of §P§. Write 
$$
P(X) = P_0(X)^\alpha R(X), \quad \mtext{with §R\in \calK[X]§, §\alpha\in\bbN^*§ 
and §R§ coprime to §P§ in §\calK[X]§.}
$$
There exists an infinite set §S\incl \calR§ such that, for every §m,m_1, m_2\in S§, 
with §m_1\ne m_2§,
\beginlist
\item {(i)} §P_0(m_1)§ is coprime to, and distinct from, §P_0(m_2)§; 
\item {(ii)} §R(m)\in \calR§;
\item {(iii)}§P_0(m)§ is coprime to §R(m)§ in §\calR§;
\item {(iv)} if §P_0§ has no multiple roots, 
§P_0(m)§ is either a unit, or is multiple of a certain prime 
§p§, but not of §p^2§; in the case where §\calU§ is finite, it
can be assumed that §P_0(m)\nin \calU§;
 \item {(v)} if §\,\content\; \calR§ is finite, §P_0(m)\in \calU§;
\item {(vi)} If §\calR = \bbZ§, §S§ can be supposed to be contained in
 §\bbZ_+§, or §\bbZ_-§.
\endlist
\endproclaim 

\proof
To satisfy to (i), (iv), (v) and (vi), apply \thm \proclref{main}
to §P_0§ in order to obtain a set §S§ satisfying the conclusion of the theorem.

We could furthermore assume that 
§P_0(m)§ is coprime to the denominators of the coefficient of §R§,
for every §m\in S§ (\cor \proclref{almostNote}).

Then §P(m)\in \calR§, §P_0(m)\in \calR§, and §P_0(m)§ is not
divisible by any prime factor of the denominators of §R§, for every §m\in S§. 
As a result, §R(m)§ must belong to §\calR§, which shows (ii).

Now, since §P_0§ and §R§ are coprime in §\calK[X]§, there exist §f,g\in \calR[X]§
such that 
$$
P_0(X) f(X) + R(X)g(X) = C,\with C\in \calR^*. \autoeqno[Bezout2]
$$ 
Eliminating from §S§ those elements §m§, in finite number,
such that §P_0(m)§ is not coprime to §C§, 
it can be assumed that §P_0(m)§ is coprime to §C§ for every §m\in S§.
From \eqref{Bezout2}, it follows that §P_0(m)§ is coprime to §R(m)§, as 
contended in (iii). 
\endproof

Assume that §\calR = \bbZ§.
It is often needed, in practice, that the elements of set~§S§ be chosen inside an 
arithmetic progression 
$$
\calA = \set{a + kb \st k\in \bbN\ (b\ne 0)}.
$$ 
For example it may be desired that all the elements of §S§ be even, odd, or multiple 
of a given number.
 
Denoting §Q(X) = a + bX§, it will be shown in section \sref{necessary_cond} that 
$$
\calP(P\circ Q) = \calP(P) \union \bigset{p\in \content\ \bbZ \st p \divis \gcd\set{P(a), b}}.
$$
Applying theorem \proclref{main} to §{1\over \gamma }P\circ Q§, 
with §\gamma=\gcd\set{P(a), b}§,
we have the following result: 

\corollary cor3
Assume that §\calR = \bbZ§.
In \thm \proclref{main}, the elements of §S§ can be chosen inside an arithmetic progression 
§\calA = \set{a + kb \st k\in \bbN\ (b\ne 0)}§, provided §\calA§ satisfy the necessary condition 
$$
\gcd\set{P(a), b} = 1.
$$
In general, if §P§ is totally primitive and if §\gamma = \gcd\set{P(a), b}§, 
there exists an infinite subset §S§ of §\calA§ such that, for every §m\in S§,  
$$P(m) = \gamma v_m$$
the numbers §v_m§ being mutually coprime.  
\endproclaim
\medbreak


The next assertion is a partial converse of 
\thm \proclref{main}, showing that if the property 
of the elements of §\set{P_i(m)}_m§ of being pairwise
coprime holds simultaneously whenever §m§ belongs to 
some set §T§,
then the conditions of the theorem are fulfilled.

\theorem secondThm
Assume that §P_1§, §P_2,\ldots§ §P_n§ are 
totally primitive polynomials in §\calK[X]§.  
If there exists an infinite set §T§ such that, for every §m_1,m_2\in T§ with 
§m_1\ne m_2§ and every §i§, §P_i(m_1)§ is coprime to §P_i(m_2)§, 
then the product of the 
polynomials §P_i§ is totally primitive, and \thm \proclref{main} holds.  
It is even possible to choose §S§ as a subset of §T§,
in order for the three first assertions of 
\thm \proclref{main} to hold.
\endproclaim 
 
\proof
Let §m\in T§,  §P = P_1(X)\cdots P_n(X)§ and §\alpha = P(m)§.  
For every §i§, only finitely many elements of §P_i(T)§ may share a prime 
factor with §\alpha§.
As a result, since §T§ is infinite, 
one can choose §m'\in S§ such that 
§P_i(m')§ is coprime to §\alpha§ for every §i§. 
Hence §P(m')§ is coprime to §P(m)§, 
showing that product of the §P_i(X)§ is totally primitive. 
This proves the first assertion of \thm \proclref{secondThm}.

In fact, we can repeat this argument to find §m''\in T§ such that 
§P(m'')§ be coprime to §\alpha' = P(m)P(m')§.

Iterating this process, an infinite set §S§
can be extracted from §T§, such that §P(m_1)§ is
coprime to §P(m_2)§ whenever §m_1, m_2\in S § and 
§m_1\ne m_2§. Thus, the two first assertions of \thm \proclref{main}
holds for §S§. 

To prove that the third assertion holds as well, 
we have to demand furthermore
that  §P(m')§, §P(m'')\ldots§ be coprime to the elements of a certain finite set 
§\calQ \incl \calR§, 
which is licit (\cor \proclref{almostNote}). 
Namely, with the notations of \thm \proclref{main}, 
for every §i,j\in \calE§ with §i\ne j§,
there exists §f_i, g_i\in \calR[X]§ such that 
$$
P_i(X) f_i(X) + P_j(X) g_i(X) = C_{i,j}, \With C_{i,j}\in \calR.
$$
We define §\calQ = \set{C_{i,j}}_{i,j}§.

So, if §m',m''\ldots m^{(k)}§ are of the aforementioned form, 
it is visible that §P_i(m^{(k)})§ 
is coprime to §P_j(m^{(k)})§, since both numbers are coprime to §C_{i,j}§.
This shows our contention.
\endproof

The following example will show how to manipulate some of the
the concepts introduced so far.

\example  
\emph{Let §Q(X)\in \bbQ[X]§; assume that for every §k\in \bbN§, §Q(k)§ is a §n§-th power
in §\bbQ§.
To show that §Q§ is a §n§-th power in §\bbQ[X]§.}
\smallbreak

If §Q§ is constant, the result is trivial. Otherwise,
Write 
$$
Q(X) = {A\over B} P(X),
$$
with §A§ coprime to §B§, and §P§ totally primitive and non constant (\prop \proclref{repr}).
Apply \cor \proclref{almostNote} to §P§ in order to obtain a set §S§ such that
§P(m)§ is coprime to §A§ and §B§, for some or other §m\in S§.
Since §Q(m)§ is a power of §n§, so must be §A§, §B§ and §P(m)§ because
these numbers are mutually coprime. 
Then §P(k) = {B\over A}Q(k)§ is also a power of §n§, for every §k\in \bbN§.

Now, let §P_0§ be a totally primitive prime factor of §P§. Write 
$$
P(X) = P_0(X)^\alpha R(X),\ \mtext{with §R\in \bbQ[X]§ and §R§ coprime to §P_0§.}
$$
Apply \cor \proclref{FactManip} in order to obtain a set §S§
and primes §p§ satisfying the conclusions of this corollary,
with respect to §P§, §P_0§ and §R§. 

Thus, for some §m\in S§, §P_0(m)§ is coprime to §R(m)§, 
§P_0(m)§ is divisible by §p§ but not by §p^2§,
§P(m)§ is an §n§-th power, and we have 
$$
P(m) = P_0(m)^\alpha R(m).
$$
It follows that §\alpha§ is multiple of~§n§.

This holds for every totally primitive factor of §P§, hence,
the primitive factorization of §P§ can be written
$$
P = {1\over \Gamma} P_0^{\alpha_0} P_1^{\alpha_1} \ldots P_n^{\alpha_n},
\quad \mtext{with §n\divis \alpha_i§ for every §i§}
$$
(\prop \proclref{unique_primitive_fact}).
But for some or other §m\in S§, §P(m)§ and §\prod_i P_i^{\alpha_i}(m)§ are both §n§-th powers; hence so is §\Gamma§.
Conclude that §P§, and hence §Q§, is a §n§-th power in §\bbQ[X]§.
\endproof


\section proof Proof of Theorem \proclref{main}

Assertion (v) in \thm \proclref{main} follows simply from the infiniteness of §S§ and
from \cor \proclref{almostNote}, applied to the (finite) set of primes of §\calR§.  

Regarding assertion~(vi),
the proof for §S \subset \bbZ_+§ or §S\subset \bbZ_-§, is the same, \emph{mutatis mutandis}, 
as for §S\incl \calR§: it suffices to exchange §\calR§ with §\bbZ_+§ (\resp §\bbZ_-§)
when necessary.  
So, we have only to prove assertions (i)--(iv) of \thm \proclref{main}.

We shall need the following lemma.

\lemma pnu2p
Assume that §P\in \calK[X]§ is an integrally valued polynomial, §p§ is a prime, 
and §a,b\in \calR§. 
For a sufficiently large exponent §\nu_p\in \bbN^*§, 
$$
{\rm if}\  a \congr b \mod {p^{\nu_p}},\  
{\rm then}\ P(a) \congr P(b) \mod p;
$$
\endproclaim

\proof
Write §P§ in the form 
$$
P = {1\over B} Q(X),
$$
with §B\in \calR§ and §Q\in \calR[X]§.
Choose §\nu_p§ larger
than the exponent of the greatest power of §p§ dividing §B§, 
and assume that §a\congr b \mod {p^{\nu_p}}§. 
Then §Q(a)\congr Q(b)\mod {p^{\nu_p}}§, \thatis
$$
Q(a) = Q(b) + kp^{\nu_p}
$$
for some §k\in \calR§; hence
$$
P(a) = {Q(a)\over B} = {Q(b) + kp^{\nu_p}\over B} = P(b) + {kp^{\nu_p}\over B}.
$$ 
Since §P(a)§ and §P(b)§ belong to §\calR§, so does the element §kp^{\nu_p}/B§, which 
is a multiple of §p§ thanks to the choice of §\nu_p§.
\endproof


Let §P_i§, §P§, §\calE§ and §\calI§ be as in the statement of \thm \proclref{main}.

For every §i,j\in \calE§ with §i\ne j§, 
there exists §f_i,g_i\in \calR[X]§
such that 
$$
P_i(X) f_i(X) + P_j(X) g_i(X) = C_{i,j}, \With C_{i,j}\in \calR^*. \autoeqno[Bezout_first]
$$

Also, there exists 
§K_i\in \calR^*§ such that §K_i P_i(X)\in \calR[X]§, and we denote
$$Q_i(X) = K_iP_i(X).$$

Now, for every §i\in \calI§,
§P_i§ and §P'_i§ have no common factor in §\calK[X]§,
hence so do §Q_i§ and §Q'_i§.
Thus, for suitable polynomials 
§\varphi_i,\psi_i\in \calR[X]§ and §D_i\in \calR^*§, there holds
$$
Q_i(X)\varphi_i(X) + Q'_i(X)\psi_i(X) = D_i. \autoeqno[Bezout_main]
$$

Let
$$
\calC  = \content \set{C_{i,j}, D_i}_{i\ne j}.
$$

Set §S_{0} = \emptyset§, and assume 
inductively that subsets
$$
S_0 \subset S_1 \cdots \subset S_k
$$
of §\calR§ have been constructed, with §\absv{S_k} = k§, such that
§S_k§ satisfies the assertions concerning §S§ in 
Theorem~\proclref{main}, except being infinite.
We denote 
$$S_k = \set{m_1, \ldots m_k}\mtext{for every §k\geq 1§, and}
\calM_k = \content\set{P(m)\st m\in S_k}. 
$$
For our induction, we assume furthermore that
$$
\calM_k \inter \calC = \emptyset. \autoeqno[indAssumption]
$$
Finally,we denote by §\Omega§ the union of the three following sets:
$$
\eqalign{
& \set{x \in \calR \st P_i(x) = P_i(m_r), \mtext{all §i§, §r\leq k§\!\!}},\cr
& \set{x \in \calR \st P(x) = P(m_r),  \mtext{§r\leq k§\!\!}}, \cr 
&\mtext{and} \set{x \in \calR \st P_i(x) = P_j(x), \mtext{§i,j\in \calE§, §i\ne j§\!\!}}.\cr
}
$$ 
These sets are finite, hence so is their union §\Omega§.

We now proceed to the construction of §S_{k+1}§, \thatis we have to built 
§m_{k+1}§.

If §\calM_k\union  \calC = \emptyset§, we choose §m§ at random 
inside §\calR\diff (S_k\union \Omega)§.
 
Otherwise §\calM_k\union \calC\ne \emptyset§.
For each §p\in \calM_k§ (if any), there exists §a_p\in \calR§ such that 
$$
P(a_p) \not\congr 0 \mod p,
$$
because §\calP(P) = \emptyset§ by hypothesis.

According to Lemma \proclref{pnu2p}, 
if §\nu_p\in \bbN^*§ is sufficiently large, there holds, for 
every §a\in \calR§, 
$$
a \congr a_p \mod {p^{\nu_p}} \implies  
P(a) \congr P(a_p) \mod p;
$$

Similarly, for each prime §q\in \calC§ (if any),
there exists §\alpha_q\in \calR§ such that
$$
P(\alpha_q) \not\congr 0 \mod q,
$$
and §\nu_q\in \bbN§ such that, 
$$
a\congr \alpha_q \mod {q^{\nu_q}} \implies P(a)\congr P(\alpha_q)\mod q.
$$  

Use the Chinese remainder theorem to find an element §m\in \calR§
congruent to each §a_p§ modulo §p^{\nu_p}§, 
and to each §\alpha_q§ modulo §q^{\nu_q}§. 
Then §P(m)§ is congruent to §P(a_p)§ modulo §p§, and to §P(\alpha_q)§ modulo §q§,
hence §P(m)§ is coprime to every element of §\calM_k \union \calC§.  

We see that whether or not §\calM_k\union \calC = \emptyset§,  
the elements 
$$P(m_1), P(m_2)\ldots P(m_k), P(m)$$ are mutually coprime, 
and 
$$
\content \set{P(m_1),\ldots P(m_k),P(m)} \inter \calC = \emptyset.
$$ 
Moreover \eq \eqref{Bezout_first} shows that §P_i(m)§ and §P_j(m)§ are coprime,
for every §i,j\in \calE§, §i\ne j§, because they are
coprime to §C_{i,j}§. 

Nevertheless, in the case where §\calM_k\union \calC\ne \emptyset§,
§P_i(m)§ may not be distinct from §P_j(m_k)§ for all §m_k\in S_k§.
To obtain distinctness and to prove the other assertions as well, we have to use the freedom provided by the Chinese
remainder theorem, and refine the choice of §m§ above.
 
Denote by §J_0§ the principal ideal generated by the product of the 
§p^{\nu_p}§ and the §q^{\nu_q}§.
Then every element §m'\in m+J_0§ preserves the above properties.
But §J_0§ is infinite, while §\Omega§ is finite.  
Therefore we can choose §m'\in (m+J_0) \diff \Omega§ such that §P_i(m')\ne P_j(m')§ if §i\ne j§,
§P(m')\ne P(m_r)§ and §P_i(m')\ne P(m_r)§ for every §r\leq k§. 
Then exchanging if necessary §m§ with §m'§, the required distinctness of §P_i(m)§ with §P_j(m)§ etc. can now be assumed.

Thus, with §m_{k+1} = m§, the above construction (\eq \eqref{Bezout_main} unnecessary) would suffice to define a set 
§S_{k+1}§ that satisfies parts (i)--(iii) of \thm \proclref{main}, related to the coprimality of the values, as well as the 
induction assumption \eqref{indAssumption}.

To prove the missing parts of the theorem, we shall use again the freedom provided by the Chinese
remainder theorem, to refine even more the element §m§.

With §J_0§ as above, observe that for every §P_i§ with §i\in \calI§, either there is a prime element 
§\pi_i§ dividing some element in §P_i(m + J_0)§, or §P_i(m+J_0)§ is a subset of §\calU§, in which case we shall
we set §\pi_i = 1§.
Notice that this last case cannot hold whenever §\calU§ is finite, 
because §m+J_0§ is infinite, while there are at most finitely many 
solutions to the equations §P_i(x)=u \in \calU§.

It should also be observed that the §\pi_i§ are coprime, according to the above part of the proof.

Denote by §\calI'§ the set of indices §i\in \calI§ for which §\pi_i^2§ divides §P_i(m)§ (hence also §Q_i(m)§), 
and by §\calI''§ the complementary of §\calI'§ in §\calI§.
 
Let §\frP § be the product of the prime elements §\pi_i§ for every §i\in \calI'§, of the elements
§\pi_i^2§ for every §i\in \calI''§, of the 
§p^{\nu_p}§  for every §p\in \calM_k§, and of the 
§q^{\nu_q}§ for every §q\in \calC§. 
Notice that §\frP§ is divisible by §\pi_i§ but not by §\pi_i^2§, for every
§i\in \calI'§ and §\pi_i\ne 1§, while §\frP§ is divisible by §\pi_i^2§, for every §i\in \calI''§.
Moreover, we have §m+\frP\in m+J_0§.

For every §i\in \calI'§, it is visible from \eq \eqref{Bezout_main}
that every prime shared by both §Q_i(m)§ and §Q'_i(m)§ must divide §D_i§.
Since §\pi_i§ is coprime to §D_i§ by construction, §\pi_i§ cannot 
divide both §Q_i(m)§ and §Q'_i(m)§.

Applying Taylor's formula up to the second order\footn{%
For a monomial §X^n§, the binome formula gives
§(X+\pi)^n - X^n = n X^{n-1}\pi \mod {\pi^2}§,
hence the formula is true in this case. 
By linearity, it is true for every 
polynomial §Q§.
}, 
we have
$$
Q_i(m + \frP) - Q_i(m) \congr \frP Q'_i(m) \mod {\frP^2}. \autoeqno[Qi2Pi]
$$
Consequently, if §i\in \calI'§ and §\pi_i\ne 1§,
§Q_i(m+\frP)§, and hence §P_i(m+\frP)§, cannot be multiple of §\pi_i^2§, since
§Q_i(m)§ is, but §\frP§ is not.

Also, if §i\in \calI''§, then §Q_i(m+\frP)§, and hence §P_i(m+\frP)§,
cannot be multiple of §\pi_i^2§, since
§\frP§ is, but §Q_i(m)§ is not.

On the other hand, §Q_i(m+\frP)§, and hence §P_i(m+\frP)§, \emph{must} be multiple of §\pi_i§
for every §i\in \calI§, since §Q_i(m)§ and §\frP§ are.  

Thus, for every §i\in \calI§, §P_i(m+\frP)§ is divisible by §\pi_i§, but not by §\pi_i^2§, 
unless §\pi_i = 1§. 

Denote by §J_1§ the principal ideal generated by the product of the 
§\pi_i^2§ for every §i\in \calI§, of the elements
§p^{\nu_p}§  for every §p\in \calM§, and of the 
§q^{\nu_q}§ for every §q\in \calC§. 
Set §m' = m + \frP§.
Let us show that every element in §m' + J_1§ 
preserves all the above properties.

First, we have §\frP + J_1\incl J_0§, hence §m' + J_1\incl m + J_0§.
This means that the coprimality properties above are preserved for 
the elements of §m'+J_1§.

Second, Taylor's formula implies that, for every §i\in \calI§,
$$
Q_i(m' + J_1) - Q_i(m') \incl J_1 \incl \pi_i^2\calR. 
$$
Therefore, for every §i\in \calI§ such that §\pi_i\ne 1§, 
the elements in §Q_i(m' + J_1)§, and hence of §P_i(m'+J_1)§, must be multiple of §\pi_i§, but not of §\pi_i^2§, 
whenever §\pi_i\ne 1§.
Thus this property is preserved as well.
Finally, if §\pi_i = 1§, since §m'+J_1\incl m + J_0§, we have §P(m'+J_1)\incl \calU§,
which ends the proof of our contention.

We can now use the infiniteness of §m'+J_1§ to find an element §m_{k+1}\in (m'+J_1)\diff \Omega§
satisfying not only the above properties, but also the required distinctness properties:
\beginlist
\item{1)} §P_i(m_{k+1})\ne P_i(m_r)§,  for every §i\in \calE§ and §r\leq k§;

\item{2)} §P(m_{k+1})\ne P(m_r)§, for every §r\leq k§;

\item{3)} §P_i(m_{k+1})\ne P_j(m_{k+1})§ for every §i,j\in \calE§ with §i\ne j§.
\endlist

The proof of \thm \proclref{main} is complete. \endproof

\comment
Regarding assertion (v), we only sketch the proof, as it needs only iterating an adjustment of the
above method.

More precisely, in place of considering a single prime element §\pi_i§ dividing §P_i(m)§, we consider 
all the prime elements §\pi_{i,j}§ dividing §P_i(m)§.
We have to replace set §\calI'§ above by the set of all indices §(i,j)§ such that §P_i(m)§ is 
multiple of §\pi_{i,j}§, but not of §\pi_{i,j}^2§, and let §\calI''§ be its complementary inside the set of all
possible indices §(i,j)§.

Now, following exactly the proof above, \emph{mutatis mutandis}, it can be deduced that 
§P_i(m+\frP)§ is divisible by §\pi_{i,j}§, but not by §\pi_{i,j}^2§ for all §(i,j)§, unless 
§P_i(m+J_0)\incl \calU§ (in what case we set §\pi_{i,j}= 1§).

If at this stage §P_i(m+\frP)§ is square-free for all §i§, we end the proof as above, 
\emph{mutatis mutandis}, by showing 
that we can choose §m_{k+1}§ inside §m'+J_1§ (with §m'=m+\frP§).  

Otherwise, §P_i(m')§ is not only multiple of §\pi_{i,j}§ but also of new prime elements
§\pi'_{i,j}§.
We can repeat the above process to obtain an element §m''\incl m+J_0§ such that §P_i(m'')§ is divisible
by §\pi_{i,j}§ and §\pi'_{i,j}§, but not by their squares. 

If §P(m'')§ is square-free, we end the proof as above. Otherwise we repeat the same process again and again
until the resulting polynomial value be square-free. 
This must be eventually the case, because the number of prime 
elements in §\calR§ is finite by hypothesis. 
\endproof
\endcomment


\section necessary_cond Note on Prime moduli annihilating a polynomial 

This section concerns integer-valued polynomials, \thatis we assume §\calR = \bbZ§.

It is well known that the integer-valued polynomials form a 
§\bbZ§-module, for which the set of binomial polynomials is a basis.
More precisely, assuming §\deg(P) = n§, Taylor's formula for finite 
differences gives
$$
P(X)= a_0 + a_1{X \binomial 1} + a_2{X \binomial 2} 
+ \cdots + a_n {X\binomial n}, 
$$
where 
$$
\eqalign{
a_0 &= P(0), \cr
a_1 &= P(1)-a_0, \cr
a_2 &= P(2) - a_0 - a_1{2 \binomial 1},\cr
&\vdots\cr
a_n &= P(n) - a_0 - a_1 {n\binomial 1} - a_2 {n\binomial 2} - \cdots 
- a_{n-1}{n\binomial n-1}.\cr
}
$$
Thus, §a_i\in \bbZ§, showing the above contention.

Thanks to this formula, it is also immediate that if §P(0)§, §P(1)\ldots P(n)§ 
share a common factor §m§, then §m§ divides the §a_i§, hence also 
the whole set of values of §P§ (the converse
is obvious). In other words, 
$$
\gcd\bigset{P(0),P(1), \ldots P(n) } = \gcd\set{a_0, \ldots a_n} = 
\gcd\, P(\bbZ).
$$
This provides a simple mean to determine the prime moduli annihilating §P§, 
\thatis §\calP(P)§.

This principle may be useful as a practical test, but is not very 
enlightening, and may be insufficient for theoretical developments.

To obtain another, somewhat more explicit, criterion, let us first restrict ourselves to 
polynomials §P\in \bbZ[X]§.
Moreover, since it is not unusual for a polynomial to be required to take 
its values inside an 
arithmetic progression (\eg the odd or even numbers), we shall consider sets 
$$
\set{P(m)}_{m\in \calA} \With 
\calA = \set{a + kb \st k\in \bbN} \ (b\ne 0), 
$$ 
and denote by §\calP(P, \calA)§ the set of prime §p§ such that §P§ vanishes identically modulo §p§ on 
§\calA§ (hence all the elements of §P(\calA)§ are multiple of §p§).
If §\calA§ is trivial (§b = 1)§, §\calP(\calA)§ is denoted, as in previous sections of this 
article, by §\calP(P)§.

Of course, choosing §a = 0§ and §b = \pm 1§ leads to the set of primes
for which §P§ vanishes identically on §\bbZ_+§ and §\bbZ_-§ modulo §p§.  

\proposition necessary_Z
Assume that §P\in \bbZ[X]§, with §\deg(P) = n§.
The set §\calP(P, A)§ is finite, and is the union (not necessarily disjoint) 
of the following three sets: 
\beginlist \itemindent=3\parindent
\item{§\calP_1(P, \calA)§:} The prime factors common to all the coefficients of §P§;
   \item{§\calP_2(P, \calA)§:} the prime factors common to §P(a)§ and §b§;
   \item{§\calP_3(P, \calA)§:} the primes §p§ dividing §P(0)§ and not larger than §n§,
    for which the reduction of §P§ modulo §X^p -X§ in §\bbZ_{/p} \,[X]§ is null. 
\endlist
\endproclaim

\proof 
It is obvious that if a prime §p§ divides all the coefficients of §P§, then it 
divides §P(m)§ for every §m§, hence §\calP_1\incl \calP(P, \calA)§.
 
If §p§ is a prime factor dividing §\gcd\set{P(a), b}§, 
then 
$$
P(a) \congr 0 \mod p \And a + kb \congr a \mod p.
$$ 
Hence 
$$
P(a+kb) \congr P(a) \congr 0 \mod p.
$$
Therefore §\calP_2 \incl \calP(P, \calA)§.

If §p§ is such that the reduction of §P§ modulo §X^p - X§
is null, then §P(k) \congr 0 \mod p§ on the whole of §\bbZ§, and \emph{a fortiori}
§p§ belongs to §S§. Hence §\calP_3\incl \calP(P, \calA)§.

Conversely, assume that §p§ is not a prime factor common to all the coefficients 
of §P§ (so §P§ is not null modulo §p§), and that §p§ does not divide §\gcd\set{P(a), b}§.  

If §P§ vanishes identically on §\calA§ modulo §p§, there must hold in particular 
§P(a) = 0 \mod p§; therefore §p§ must be coprime to §b§ according to our assumptions. 
As a consequence, §b§ is invertible modulo §p§, and §a+kb§ takes all the 
values §0, 1, 2, ... p-1§ modulo §p§ as §k§ varies in §\bbN§. 
Hence §P(k) = 0 \mod p§ for every integer §k§, which is possible only if the 
reduction of §P§ modulo §X^p - X§ is null in §\bbZ_{/p}\, [X]§. 

Furthermore, since §P§ is not null modulo p, it cannot have more roots in
§\bbZ_{/ p}§ than its degree §n§, that is, §p\leq n§.

Finally, since the reduction of §P§ modulo §X^p-X§ in §Z[X]§ leaves §P(0)§ unchanged, 
§P(0)§ must be congruent to §0§ modulo §p§.

This shows that §\calP(P, \calA) \incl \calP_1\union \calP_2\union \calP_3§.  
And since §\calP_1§, §\calP_2§ and §\calP_3§ are obviously finite, so is §\calP(P, \calA)§.
\endproof

In general, assume that §P§ is an integer valued polynomial.
If §P§ is written in ``binomial form'', then §\gcd\, P(\bbZ)§ 
is equal to the §\gcd§ of the coefficients of §P§.

Otherwise, §P§ may be written in the form
$$
P = {1\over \gamma} Q, \With Q\in \bbZ[X] \And 
\gamma\in \bbZ.
$$
Obviously, there must hold
$$
\content \set{\gamma} \incl \calP(Q, \calA),
$$
hence 
$$
\calP(P, \calA) \incl \calP(Q, \calA) = \calP_1(Q, \calA) \union \calP_2(Q, \calA) 
\union \calP_3(Q, \calA),
$$
with the above notations.

But a criterion that would provide the exact set §\calP(P, \calA)§ would be more complex and 
probably cumbersome. This question is related to the determination of the integer valued 
polynomials which 
vanish identically modulo a composite number §m§. 
See Kempner~(\cite{Kem}) and Singmaster~(\cite{Sin}) for more information on this topic. 

\example
\emph{To show that §P(X) = X^6 + 3X^5 + X^2 + X + 21§ is totally primitive.}
 
Apply the previous theorem to §P(X)§, 
with §\calA = \bbZ_+§, the trivial arithmetic progression (§a=0§, §b = 1§).  
We have §\calP_1(P, \calA) = \emptyset§ since §P§ is primitive in the usual sense,
and of course, §\calP_2(P, \calA) = \emptyset§. On the other hand, 
the only prime factor not larger than 
§\deg(P)§ that divides §P(0)§ is §3§, and the reduction of §P(X)§ modulo §X^3 - X§
in §\bbZ_{/3}§ is §2X^2 + X \ne 0§.
Thus §\calP_3(P, \calA) = \emptyset§.  


\section conjectures Bunyakovsky and Schinzel conjectures restated and extended

We have set a convenient framework to deal with questions related to 
the values of polynomials. 
It may be interesting to restate and extend some famous conjectures
in this context.

Assume that §P§ is a totally primitive polynomial, whose value set contains 
infinitely many prime numbers. 
Obviously, §P§ must be non constant, with positive leading coefficient.
Denote by §S\incl \bbZ§ an infinite set such that §P(m)§ is prime for every §m\in S§.

If §P_0§ is a totally primitive prime factor of §P§, write 
$$
P(X) = {1\over \Gamma}P_0(X)R(X),\quad \mtext{with §R\in \bbZ[X]§ and §\Gamma\in \bbN^*§.} 
$$
Then for every §m\in S§,
$$ P(m) = {1\over \Gamma} P_0(m) R(m),$$
hence §P(m)§ must divide §P_0(m)§ or §R(m)§, since §P(m)§ is prime.
Furthermore, for infinitely many §m\in S§, §P(m)§ is coprime to §\Gamma§,
therefore §P(m)§ divides either §P_0(m)§ or §R(m)§, but not both, and the one
that is not divisible by §P(m)§ must divide §\Gamma§. 
Thus, for infinitely many §m\in S§, §P_0(m)§, or §R(m)§, divide §\Gamma§.
Since the number of divisors of §\Gamma§ is finite, this is possible only if §P_0§ or 
§R§ is constant. But §P_0§ has been supposed to be a totally primitive prime factor of §P§, 
hence §P_0§ is not constant, and §R§ must be constant.
This shows that §P§ \emph{must} be irreducible in §\bbQ[X]§.

Thus, in our framework, Bunyakovsky's conjecture can be extended as follows:

\conjecture Bunyakovsky 
If §P\in \bbQ[X]§ is a non constant and totally primitive polynomial, of positive leading coefficient, 
then the value set of §P§ contains infinitely many prime numbers.  
\endproclaim

This conjecture is a consequence of the more general ``extended'' Schinzel conjecture
we shall now deal with.

Let §P_1, P_2,\ldots P_n§ be totally primitive polynomials. 
Let us find necessary conditions for the existence of an infinite set §S§ such that 
the §P_i(m)§ be simultaneously prime for every §m\in S§. 

Assuming this assertion holds, we first observe
that §P_i(S)§ must be infinite for every §i§, 
since §P_i(x) = p§ has only finitely many solutions for each §p§.
This implies, as we have just seen, that the §P_i§ must be non constant,
irreducible in §\bbQ[X]§, 
and have positive leading coefficients.  

Now, since the number of solution of §P(x) = p§ is finite for every §p§, 
we can eliminate those §m\in S§ that give rise repeated values of §P_1(m)§, while keeping 
§S§ infinite,  
and hence assume that §S§ is such that the §P_1(m)§ are distinct, whenever §m§ varies in §S§.
Similarly, we can eliminate those §m§ in this new set §S§ that give rise repeated values of §P_2(m)§, 
and continuing so, we finally get an infinite set §S§ such that §P_i(m_1)\ne P_i(m_2)§ for every §i§,
whenever §m_1\ne m_2§, with §m_1, m_2\in S§.

According to \thm \proclref{secondThm}, it follows that the product 
§P_1P_2\cdots P_n§ is totally primitive.

These are the desired necessary conditions, which make plausible the following conjecture, 
very similar to Schinzel's hypothesis~H: 
\medbreak

\conjecture Schinzel 
Assume that §P_1, P_2,\ldots P_n\in \bbQ[X]§ are totally primitive and non constant polynomials, 
irreducible in §\bbQ[X]§, with positive leading coefficients. 
If the product of the §P_i§ is totally primitive, then there exists an infinite set §S§ such 
that § P_i(m)§ is a prime number, for every §i§, and every §m\in S§.
\endproclaim 

From the remarks above, \thm~\proclref{secondThm}, and the fact that powers 
of totally primitive polynomials are totally primitive,
the following enlightening statement would follows from the above conjecture. 

\corollary[(strong form)] Schinzel2
Suppose that §P\in \bbQ[X]§ is a non constant, totally primitive polynomial, of positive leading coefficient.
Assume that §P§ admits a factorization into a product of powers of totally primitive distinct 
prime factors: 
$$
P = P_1(X)^{\alpha_1} P_2(X)^{\alpha_2}\ldots P_s(X)^ {\alpha_s}.
$$

Then there exist an infinite set §S\incl \bbN§ such that, for every §m\in S§, 
$$
P(m) = p_1^{\alpha_1} p_2^{\alpha_2}\ldots p_s^{\alpha_s},\quad  
\with p_i\in \bbP,\  p_i\ne p_j\ \forall i\ne j,
$$
and such that §P(m_1)§ is coprime to §P(m_2)§, for every §m_1, m_2\in S§.
\endproclaim


\section PrincConj    Related conjectures

It is natural to suspect that the Bunyakovsky and Schinzel conjectures, with obvious 
adaptations, hold in principal ideal domains.\footn{See for example the 
work of Bender et all. (\cite{BW}) for an analogue of the Schinzel 
hypothesis for polynomials over §\bbF_q[t]§.}
Unfortunately, this is not the case in general; as shown in \cite{CG},
there exist totally primitive and irreducible polynomials over §\bbF_q[t]§, whose
all values are reducible in §\bbF_q[t]§.

For example, assume that $\bbF_p$ is the field with $p$ elements, $p$ odd. 
Choose an integer $a$ with $1< a < 4p$, such that $a$ be coprime to $p(p-1)$ 
(for example, $a = 2p - 1$ is suitable). Let $t$ be transcendental over $\bbF_p$,
and set $\calR = \bbF_p[t]$. Define $$P(X) = X^{4p} +t^a\in \calR[X],$$
which can be shown to be totally primitive, and irreducible (see \cite{Lang},
\thm.~9.1).
Then for every $f\in \bbF_p[t]$, $P(f)$ is reducible in $\bbF_p[t]$.

Here is a direct derivation of this example, 
which was pointed out by Prof.~Peter Mueller.{\makechar\#\footn{See mathoverflow, 
\beginhyperlink{https://mathoverflow.net/questions/412271/%
proving-that-polynomials-belonging-to-a-certain-family-are-reducible?noredirect=1#comment1056899_412271}%
Proving that polynomials belonging to a certain family are reducible.\endhyper}} 
We shall use a result of Swan (\cite{Sw}).
If $t$ divides $f(t)$ there is nothing to do. 
So assume this is not he case. Then with $F(t)=P(f(t))=f(t)^{4p}+t^a$, 
we have $F'(t)=at^{a-1}$, and $F(t)$ is separable. 
From this, one easily gets that the discriminant of $F(t)$ is a square in $\bbF_p$ 
(see the alternative formula for the discriminant at the beginning of the article of Swan). 
By Swan's result (Corollary~1), the number of irreducible factors of $F(t)$ is congruent to 
$\deg F(t)=4p\deg f$ mod $2$, hence is even, and therefore $>1$.\endproof
  
Despite the author is unaware of any counter-example in characteristic 0, it would be 
hazardous, without solid reasons, to extend the aforementioned conjectures to 
principal ideal domains of zero characteristic.

Actually, we believe these conjectures are made up of two problems, of very different nature:
stated in the more general context of principal ideal domains, the first problem is to bound the number 
of prime factors of the values of a polynomial:  
\medbreak

{\sc Problem 1:}
{\it 
Assume that §P\in \calK[X]§ is a non constant, totally primitive polynomial.
Let §\calM_n§ denotes, for every §n\in \bbN§, the (possibly empty)
set of elements §m\in \calR§ such that §\bigabsv{\content\set{P(m)}} = n§.
To find
$$
L = \inf \set{n \st \absv{\calM_n} = \infty} \in \bbN\union \set{\infty}.  
$$}%

The Bunyakovsky (extended) conjecture implies that if §\calR = \bbZ§ and §P§ is irreducible,
then §L = 1§.
The Schinzel (extended) hypothesis implies that if §\calR = \bbZ§ and §P§ has a 
totally primitive factorization, free of constant divisor, into powers of §n§ 
irreducible factors, then §L = n§.
This problem has been known for long in one or other form, and is, of course, very difficult.
\medbreak

The second problem is to ``free'' the values of a polynomial from the squares of which they are multiple. 
\medbreak

{\sc Problem 2:}
{\it
Assume §P\in \calK[X]§ is a non constant and totally primitive polynomial, whose totally primitive factorization 
is free of constant divisor, and whose totally primitive irreducible factors §P_i§ are separable.
To determine whether there exists an infinite set §S§ such that, for every §m\in S§ and every §i§, 
§P_i(m)§ is square-free in §\calR§.  
}
\medbreak

As we have seen above, the extended Schinzel conjectures provides an affirmative answer whenever 
§\calR = \bbZ§.
Also, according to \thm \proclref{main}~(v), the answer is positive whenever §\content \,\calR§ is finite.
Actually,
it seems unlikely that the infinitude of §\content\, \calR§ would ``block'' the possibility of a solution to Problem~2, 
and \thm \proclref{main}~(iv) seems to corroborate this opinion.
This is why we believe the following conjectures are true in any principal ideal domain~§\calR§.

\conjecture PrincipalConj 
Assertion (iv) in  \thm \proclref{main} can be replaced by ``§P_i(m)§ is square free''. 
\endproclaim

If we reduce the above assertion to a single polynomial, and drop the coprimality condition,
we obtain the following weaker conjecture.

\conjecture PrincipalConj2
If §P\in \calK[X]§ is a non constant, totally primitive and separable polynomial, 
the value set of §P§ contains infinitely many square-free elements in~§\calR§.
\endproclaim

We believe this problem is easiest than the first problem, albeit difficult.
The methods that need to be developed in order to ``fight'' the infinitude of primes might 
be far reaching, and an important step toward the proofs of the  Bunyakovsky and Schinzel
conjectures.

\def\gobbleit#1{}%

\bigbreak\bigskip\centerline{\bf Bibliography}
\medskip \nobreak
\frenchspacing
\tenrom
\baselineskip = 12\scpt

\bibliographystyle{plain}

\bibliography{coprimePoly}

\begin{thebibliography}{1}

\bibitem{BW}
O.~Wittenberg A.~O.~Bender.
\newblock {A} potential analogue of the schinzel hypothesis for polynomials
  with coefficients in {$\bbF_q[t]$}. {S}ubmitted.

\bibitem{CG}
R.~Gross B.~Conrad, K.~Conrad.
\newblock {I}rreducible specialzation in genus 0. {S}ubmitted.

\bibitem{Kem}
A.~J. Kempner.
\newblock {P}olynomials and {T}heir {R}esidue {S}ystems.
\newblock {\em Transactions of the American Mathematical Society},
  22(3):267--88, 1921.

\bibitem{Lang}
S.~Lang.
\newblock {\em {A}lgebra}.
\newblock New York, 1970.

\bibitem{Sin}
D.~Singmaster.
\newblock {O}n {P}olynomial {F}unctions (mod m).
\newblock {\em Journal of Number Theory}, 6(5):345--352, 1974.

\bibitem{Su}
B.~Sury.
\newblock {P}olynomials with {I}nteger {V}alues.
\newblock {\em Resonance}, 6(9):46--60, 2001.

\bibitem{Sw}
R.~G. Swan.
\newblock {Factorization of polynomials over finite fields.}
\newblock {\em Pacific Journal of Mathematics}, 12(3):1099 -- 1106, 1962.

\end{thebibliography}


\end{document}